\documentclass[10pt]{article}
\usepackage{amsmath,amsfonts,amssymb,latexsym}

\catcode`@=11
\renewcommand\section{\@startsection {section}{1}{\z@}%
                                   {-3.5ex \@plus -1ex \@minus -.2ex}%
                                   {2.3ex \@plus.2ex}%
                                   {\normalfont\normalsize\bfseries}}
\renewcommand\subsection{\@startsection{subsection}{2}{\z@}%
                                     {-3.25ex\@plus -1ex \@minus -.2ex}%
                                     {1.5ex \@plus .2ex}%
                                     {\normalfont\normalsize\bfseries}}
\def\@seccntformat#1{\csname the#1\endcsname.\enspace}
\def\@begintheorem#1#2{\trivlist
   {\item[]\bfseries #1\ #2.}\unskip\itshape\enspace}
\def\@opargbegintheorem#1#2#3{\trivlist
     {\item[]\bfseries #1\ #2\ (#3).}
      \unskip\itshape\enspace}
\catcode`@=12
\newdimen\mydimen
\def\rightpicture#1{\par\noindent{%
        \setbox0\vtop to 0pt{\hbox{\quad#1}}%
        \global\mydimen=\wd0
        \hbox to 0pt{\hbox to \textwidth{\hss
        \vbox to 0pt{\vskip-1ex\box0}}\hss}%
        }\hangindent=-\mydimen}
\newtheorem{theorem}{Theorem}[section]
\newtheorem{corollary}[theorem]{Corollary}
\newtheorem{lemma}[theorem]{Lemma}
\newtheorem{proposition}[theorem]{Proposition}

\newenvironment{definition}
{\trivlist
\item[]\noindent{\bf Definition.}\enspace\ignorespaces}
{\endtrivlist}

\newenvironment{example}
{\trivlist
\item[]\noindent{\bf Example.}\enspace\ignorespaces}
{\endtrivlist}

\newenvironment{remark}
{\trivlist
\item[]\noindent{\bf Remark.}\enspace\ignorespaces}
{\endtrivlist}

\newenvironment{claim}
{\trivlist
\item[]\noindent{\bf Claim.}\bgroup\it\enspace\ignorespaces}
{\egroup\endtrivlist}

\def\prooftext{}
\newenvironment{proof}
{\trivlist
\item[]\noindent{\it Proof \prooftext\unskip.}\enspace\ignorespaces}
{\qquad$\Box$\endtrivlist\gdef\prooftext{}}

\def\supp{\mathop{\rm supp}}

\def\closure{\mathop{\rm Cl}}
\def\codim{\mathop{\rm codim}}

\begin{document}
\pagestyle{empty}
%
%

\null\vskip 3.2 cm
\par\noindent
%
%
{\Large\textbf{Classification of two-orbit varieties}}

\vskip 0.7 cm
\par\noindent
%
%
St\'ephanie Cupit--Foutou

\vskip 1.25 cm
\par\noindent
%
%

\bigbreak
\noindent
%
%
\section{Introduction}

Our base field is the field $\mathbb C$ of complex numbers. 
We study normal complete algebraic varieties $X$
endowed with an action of 
a connected reductive algebraic group $G$
admitting a dense orbit $\Omega$.
The stabilizer $H$ of a point of the dense orbit
will be called \emph{generic stabilizer}.
Clearly, $H$ is parabolic if and only if $X=\Omega$.

If $H$ is connected, \textsc{Borel} showed in~\cite{Bor} 
that the homogeneous space
$G/H$ has at most two ends, \emph{i.e.} the complement
$X\setminus\Omega$ has at most two connected components.
When $X\setminus\Omega$ is disconnected,
the generic stabilizer contains a maximal unipotent subgroup $U$ of $G$.
This property results from 
the connectedness (proved in~\cite{BB,H}) 
of the $U$-fixed point subset of $X$.
In this situation, the varieties $X$ share other peculiarities
which were clarified by \textsc{Ahiezer} in~\cite{A2} where he
studied and classified these objects.

This article deals with the left situation -- 
$X\setminus\Omega$ is connected --,
as we assume this
complementary to be $G$-homogeneous.
Actually, this case appears as the most natural and ``simplest'' one
to consider after the case of projective homogeneous varieties (Grassmaniann, flag
varieties, etc.).
We call naturally such varieties $X$ \emph{two-orbit varieties}.
These varieties have been intensively studied before. 
\textsc{Ahiezer} in~\cite{A1} gave a classification of the pairs $(G,H)$ 
such that $\Omega=G/H$ admits a compactification $X$ 
by one homogeneous divisor. 
This list was independently obtained by
\textsc{Huckleberry} and \textsc{Snow} \cite{HS} in
the more general context of k\"ahlerian varieties.
Later, \textsc{Brion} gave in \cite{Br1,Br2} 
a purely algebraic approach to these results
using the general theory of embeddings 
developed by~\text{Luna} and \textsc{Vust} (\emph{see}\/~\cite{LV}). 
Finally, \textsc{Feldm\"uller}
classified in~\cite{Feld} all pairs $(G,H)$ giving raise to
two-orbit varieties whose closed orbit is of codimension
$2$.

The aim of this paper is to give the complete classification of all
two-orbit varieties and to prove Luna's conjecture, namely that
all two-orbit varieties are spherical, \emph{i.e.} admit
a dense orbit of a Borel subgroup $B$ of $G$.
In particular, spherical varieties have only finitely many $B$-orbits.
We refer the reader to \cite{LV,BLV}
for an introduction to this subject.

To obtain the classification, we use two main steps.
First of all, we show that defining a two-orbit variety as being normal
is not a major restriction since the normalization map
for a complete variety with two orbits is bijective  
(\emph{see}~section~\ref{notation}).
Secondly, we give only an explicit list of all \emph{cuspidal}
two-orbit varieties from which we obtain all the two-orbit varieties
by a parabolic induction procedure explicitly described 
in section~\ref{notation}.

\bigbreak

\noindent
\textbf{Acknowledgements.}
\enspace
I would like to thank my
adviser Peter~\textsc{Littelmann} for his constant support.
I am also grateful to Michel~\textsc{Brion} for his advice and remarks on
my PhD's work in which I obtained the results presented in this paper
(\emph{see}~\cite{Cu}).

\bigbreak

\noindent
\textbf{Note added in proof.}
Alexander Smirnov found similar results while I was submitting this article. 
\enspace

\section{Notation and main results}
\label{notation}
Let $G$ be a connected reductive algebraic group.
We fix a Borel subgroup $B$ and a maximal torus $T$ in $B$. 
Let $\Phi$ be the root system of $G$,
we denote $\Phi^+\subset \Phi$ the set of positive roots relative to $B$
and $\Delta\subset \Phi^+$ the set of simple roots corresponding to $B$.
In case $G$ is simple, we enumerate the simple roots 
as in~\cite{BBKI}
and when it is more convenient, 
we use the short notation of [loc. cit.] 
for the description of positive roots.
For example, if $G$ is of type $F_4$, 
$1111$ stands for the root $\alpha_1+\alpha_2+\alpha_3+\alpha_4$. 

The subset $\langle\gamma_1$,\ldots,$\gamma_r\rangle$ of $\Phi$
will refer to the set of roots obtained 
as $\mathbb Q$-linear combinations of 
the given roots $\gamma_1,\ldots ,\gamma_r$. 

For $\alpha$ a root of $G$, we have the corresponding root space 
$\mathfrak g_\alpha\subset \mathfrak g=\mathop{\rm Lie} G$.
We let $Y_\alpha$ be any element of $\mathfrak g_\alpha\setminus\{0\}$,
then $\mathbb C Y_\alpha=\mathfrak g_\alpha$.
The 0-th root space of $\mathfrak g$ is $\mathop {\rm Lie} T=\mathfrak t$.

Throughout this article, $X$ will denote
a connected normal complete algebraic variety 
endowed with an action of $G$ such that $X$ has two $G$-orbits.
Such a variety is naturally called a \emph{two-orbit $G$-variety}.
If $G\cdot x$ corresponds to the dense orbit and
$G\cdot y$ to the closed one, we have
$X=G\cdot x\cup G\cdot y=\closure(G\cdot x)$ 
(the closure of $G\cdot x$). 
In other words, a two-orbit $G$-variety can be regarded as
a complete two-orbit embedding of the homogeneous space $G/H$ where
$H$ is a generic stabilizer.

Actually, a two-orbit variety is even projective. 
For, according to Sumihiro (\emph{see}\/~\cite{Sumi}), 
the closed orbit admits a $G$-stable
quasi-projective neighbourhood -- which has to be all of $X$.
As a weak converse, we state:
a connected projective $G$-variety with two $G$-orbits 
(say $X'=G\cdot x'\cup G\cdot y'$)
is bijective to its normalization, 
which is, in particular, a two-orbit $G$-variety.
Let us sketch the proof.
Consider the $G$-equivariant morphism 
$\pi: \widetilde X'\rightarrow X'$ 
given by the normalization of $X'$.
First of all, $\pi$ being birational, 
$\widetilde X'$ has a dense $G$-orbit isomorphic to $G\cdot x'$.
Secondly, due to the finiteness of $\pi$ together with
the connectedness of a parabolic subgroup,
we get that $\pi^{-1}(G\cdot y')$ consits of a finite union
of $G$-orbits all bijective to $G\cdot y'$.
Therefore, all these orbits are projective.
Finally, we prove, by contradiction, that actually  
$\pi^{-1}(G\cdot y')$ is $G$-homogeneous. 
Suppose that $\pi^{-1}(G\cdot y')$ consists of more than 
one (projective) $G$-orbit. By ~\cite{Bor}, 
$\widetilde X'\setminus G\cdot x'$ 
will have two connected components.
Thus,
(\emph\/{see}\/ the corresponding argument given in the introduction), 
the generic stabilizer of a point of $\widetilde X'$, 
and hence of $X'$, will contain a maximal unipotent subgroup of $G$,
say $U$ subgroup of $B$.
Consequently, we will end up with two $B$-fixed points in $X'$ 
-- a situation which can not occur.

Similarly, by considering the $G$-equivariant finite morphism
given by the projection $G/H\rightarrow G/H^\circ$,
one can prove that:
if $G/H$ has a complete two-orbit embedding, then so does $G/H^\circ$.
The converse is obviously true.
More precisely, if $Z$ is a complete two-orbit embedding of $G/H^\circ$,
$ G/H$ inherits a complete two-orbit embedding given 
as the quotient of $Z$ by $H/H^\circ$.
As a consequence, we will restrict
our study to two-orbit varieties with connected generic stabilizers, 
without loss of generality.

Before stating the main results on the two-orbit varieties,
we need to make one more remark.
Blowing down a given two-orbit variety along its closed
orbit may produce a new two-orbit variety.
We will illustrate and clarify this natural geometrical
construction via the next example.

Let $G$ be a simple group of rank~$2$ 
and $\lambda$ a dominant weight such that
$(\lambda,\alpha_1^\vee)>1$ and
$(\lambda,\alpha_2^\vee)=0$,
$(\cdot,\cdot)$ being the Killing form of $G$ and 
$\alpha_i^\vee$ the coroot associated to $\alpha_i$.
We consider the irreducible $G$-module $V(\lambda)$
associated to $\lambda$
and the weight vector $v_{\lambda-\alpha_1}\in V(\lambda)$ 
of weight $\lambda-\alpha_1$.
Take $P\supset B$ the parabolic subgroup associated to $\alpha_1$.
Then
$\closure\bigl(P\cdot [v_{\lambda-\alpha_1}]\bigr) \subset
\mathbb P\bigl(V(\lambda)\bigr)$ 
is a $P$-variety with two $P$-orbits
and 
$\closure\bigl(G\cdot [v_{\lambda-\alpha_1}]\bigr)=
G\cdot\closure\bigl(P\cdot [v_{\lambda-\alpha_1}]\bigr)$ 
is a two-orbit $G$-variety.
We have an obvious $G$-equivariant morphism 
$G\times_P \closure\bigl(P\cdot[v_{\lambda-\alpha_1}]\bigr)
\rightarrow
\closure\bigl(G\cdot [v_{\lambda-\alpha_1}]\bigr)$.
This morphism is only isomorphic on the open $G$-orbits,
the morphism on the closed $G$-orbits being
$G/B\rightarrow G/Q$
where $Q\supset B$ is the parabolic subgroup associated to $\alpha_2$.

This example is very instructive in the sense that it gives the
procedure to obtain all two-orbit varieties 
and motivates the following notion inspired by Luna's work \cite{Luna}.
In the statement of this definition (and in the rest of the text),
the terminology \emph{two-orbit $P$-variety}
for $P$ a parabolic subgroup of $§G$ will refer, by abuse of language, to a complete normal $P$-variety,
with two $P$-orbits, on which the unipotent radical of $P$ acts trivially.

\begin{definition}
A two-orbit $G$-variety $X$ is
\emph{obtained by parabolic induction from a pair} $(P,Y)$
if $P\subsetneq G$ is a parabolic subgroup and $Y$ a
two-orbit $P$-variety such that 
\smallbreak 
{\rm(i)}\enspace the radical of $P$ acts trivially on $Y$;
\smallbreak 
{\rm (ii)}\enspace there exists a
$P$-equivariant injective morphism $\varphi: Y\rightarrow X$ 
inducing a birational morphism $G\times_P
Y\rightarrow X$.
\smallbreak
A two-orbit variety is \emph{cuspidal} if
it can not be obtained by parabolic induction.
\label{definduc}
\end{definition}

\begin{proposition}
\smallbreak
{\rm (i)}\enspace The geometry of a two-orbit variety obtained
by parabolic induction from a pair $(P,Y)$
is completely determined by its closed orbit
and the geometry of $Y$. 
\smallbreak
{\rm (ii)}\enspace Each (non cuspidal) two-orbit variety 
is obtained by parabolic induction from a unique pair $(P,Z)$ such that 
$Z$ is cuspidal.
\end{proposition}

Once this result is proved (\emph{see}\/~section~\ref{induction}),
we can concentrate only on cuspidal two-orbit varieties.
We get their description as follows.
We start to show which groups can occur as generic stabilizers 
(\emph{see}\/~sections~\ref{typeI} and \ref{typeII}).
This is done mainly with combinatorial methods.
Thus, afterwards, it is just a matter of computations,
using Luna-Vust Theory (\emph{see}\/~\cite{LV}), 
to describe the corresponding complete two-orbit embeddings.
We list the pairs in this section but the embeddings only in the appendix 
to avoid too much notation in this part of the text.
Recall that we consider only the case of connected generic stabilizers
$H$, 
which is not restricted at all as 
explained in the beginning of this section.
 
\begin{theorem}

The cuspidal two-orbit varieties are obtained as
embeddings of the homogeneous spaces $G/H$ with $(G,H)$ 
in Table~$1$ and Table~$2$.
In particular, $G$ is simple or equal to $SL_2\times SL_2$, in this situation.
Moreover, each such $G/H$ has exactly one complete two-orbit embedding
described in the appendix.
\end{theorem}

\begin{corollary}
Two-orbit varieties are spherical.
\end{corollary}

\begin{proof}
As the property of being spherical is not destroyed
by parabolic induction, it suffices to check that
the homogeneous spaces $G/H$ obtained from Table 1 and 2 are indeed spherical.
It appears that these embeddings either have been already classified
by~\cite{A1,HS,Br2,Feld}
or are wonderful (thus spherical by~\cite{Luna})
of rank $2$ (\emph{see}\/~\cite{Wass}).
\end{proof} 

\begin{figure}
\begin{center}
\textbf{Table~1.\enspace Pairs of type I}
\par\nobreak\medskip\nobreak
$
\begin{array}{ccll}
A_n
&\triangleright&
\mathfrak{gl}_n
&
\\
\\
B_n
&\triangleright&
\mathfrak{so}_{2n}
&
\\
\\
B_n
&\triangleright&
        
        \mathfrak t
        \oplus
        \smash{\bigoplus_{\alpha\in\Psi}}
        \mathfrak g_\alpha ,
&\quad\hbox{with}\ 
        \Psi
        =
        \langle
        \pm \alpha_1,\dots,\pm \alpha_{n-1},\alpha_{n-1}+2\alpha_n
        \rangle
\\
\\

\vcenter{\setbox0\hbox{$B_n$}
        \copy0\medskip\hbox to \wd0{\hss or $C_n$\hss}}
&\triangleright& 
        \mathfrak t
        \oplus
        \smash{\bigoplus_{\alpha\in\Psi}}
        \mathfrak g_\alpha ,
&\quad\hbox{with}\ 
        \Psi
        =
        \langle
        \pm \alpha_2,\dots,\pm \alpha_{n-1},
        \alpha_n,\varepsilon_1+\varepsilon_n
        \rangle
\\ 
\\
C_n
&\triangleright&
        \mathfrak{sp}_2\times \mathfrak{sp}_{2n-2} 
&
\\ 
\\
C_n
&\triangleright& 
        \mathfrak t
        \oplus
        \smash{\bigoplus_{\alpha \in \Psi }}
        \mathfrak g_\alpha ,
&\quad\hbox{with}\ 
        \Psi
        =
        \langle
        \pm \alpha_2,\dots ,\pm \alpha_n,2\varepsilon_1
        \rangle
\\
\\

F_{4}
&\triangleright&
        \mathfrak{spin}_9

\\
\\
F_4
&\triangleright&
        \mathfrak t
        \oplus
        \bigoplus_{\alpha\in\Psi}
        \mathfrak g_\alpha ,
&\quad\hbox{with}\ 
        \Psi
        =
        \langle\pm\alpha_2,\alpha_3,\pm \alpha_4,
        1120\rangle
\\
\\
G_2 &\triangleright&
        \mathfrak t
        \oplus
        \smash{\bigoplus_{\alpha\in\Psi}}
        \mathfrak g_\alpha ,
&\quad\hbox{with}\ 
        \Psi
        =
        \{\alpha_1,2\alpha_1+\alpha_2,
        3\alpha_1+\alpha_2,3\alpha_1+2\alpha_2\}
\\
\\
G_2
&\triangleright&
        \mathfrak t
        \oplus
        \bigoplus_{\alpha\in\Psi}
        \mathfrak g_\alpha ,
&\quad\hbox{with}\ 
        \Psi
        =
        \{\alpha_2,2\alpha_1+\alpha_2,
        3\alpha_1+\alpha_2,3\alpha_1+2\alpha_2\}
\\
\\
G_2
&\triangleright&
        \mathfrak t
        \oplus
        \bigoplus_{\alpha\in\Psi}
        \mathfrak g_\alpha ,
&\quad\hbox{with}\ 
        \Psi
        =
        \{\pm\alpha_2,2\alpha_1+\alpha_2,
        3\alpha_1+\alpha_2,3\alpha_1+2\alpha_2\}
\\
\\
G_2
&\triangleright&
        \mathfrak{sl}_3 
\end{array}
$
\end{center}
\vfill
\begin{center}
\textbf{Table~2.\enspace Pairs of type II}
\par\nobreak\medskip\nobreak
$
\begin{array}{ccll}
A_1\times A_1 &\triangleright& \mathfrak{sl}_2 
\quad\hbox{(diagonally embedded)} 
\\
\\
\vcenter{\setbox0\hbox{$A_2$, $B_2$}
        \copy0\medskip\hbox to \wd0{\hss or $G_2$\hss}}
        &\triangleright&
        \displaystyle 
        \mathop{\rm Lie}\bigl(\ker(\alpha_1-\alpha_2)\bigr)
        \oplus 
        \mathbb{C}(Y_{\alpha_1}+Y_{\alpha_2})
        \oplus
        \smash{\bigoplus_{\alpha\in\Phi^+\setminus\{\alpha_1,\alpha_2\}}}
                \mathfrak g_\alpha
\\
\\
B_3
&\triangleright& \mathfrak g_2 
&
\\ 
\\
C_3
&\triangleright&
\vcenter{\setbox0\hbox{
$\displaystyle
        \mathfrak{sl_2}
        \oplus 
        \mathbb C
        \oplus
        \mathbb C (Y_{-\alpha_2}+Y_{\varepsilon_1+\varepsilon_2})
        \oplus
        \bigoplus_{\alpha\in \Phi^+\setminus\{\alpha_2, 
                                       \varepsilon_1+\varepsilon_2,
                            \varepsilon_1-\varepsilon_3,2\varepsilon_2\}}
           \mathfrak g_\alpha $}
\copy0\medskip\hbox to \wd0{with\quad 
$        
        \mathfrak t\cap\mathfrak {sl_2}
        \oplus
        \mathbb C
        =\mathop{\rm Lie}\bigl
        (\ker(\varepsilon_1+2\varepsilon_2-\varepsilon_3)\bigr) $
        
\hss}}
\\
\\
D_n
&\triangleright&
        \mathfrak{so}_{2n-1} 
\\ 
\\
F_4
&\triangleright&
\vcenter{\setbox0\hbox{        
$\displaystyle
        \mathfrak{g_2}
        \oplus 
        \mathbb C
        \oplus
        \mathbb C (Y_{-\alpha_4}+Y_{1232})
        \oplus
        \bigoplus_{\alpha\in\langle 1121,1221,0121,\alpha_3\rangle}
                \mathfrak g_\alpha $}
\copy0\medskip\hbox to \wd0{with\quad
$ 
        \mathfrak t\cap\mathfrak {g_2}\oplus\mathbb C
        =\mathop{\rm Lie}\bigl(\ker(\alpha_4+1232)
\bigr)$ 
\hss}}

\end{array}
$
\end{center}
\end{figure}

\section{General properties}
\label{generalities}

First of all, note that we can assume
the group $G$, acting on $X$, to be semisimple.
Indeed, the radical $R(G)$ of $G$ acts trivially on $X$.
For, if there exists a generic element of $X$ not fixed by $R(G)$,
we will obtain at least two $R(G)$-fixed points 
in the closure of its $R(G)$-orbit.
But it is a well-know fact, that these fixed points belong 
to different connected components of the $R(G)$-fixed point set in $X$ --
which is connected in our situation 
since it is equal to the closed $G$-orbit of $X$, according to

\begin{lemma}
Let $H$ be a subgroup of $G$ and $T'$ be a
torus of $G$. If $G'$ is the identity-component of the
centralizer in $G$ of $T'$, then the connected components
of $(G/H)^{T'}$, the $T'$-fixed points of $G/H$, are
exactly given by its $G'$-orbits.
\label{compos}
\end{lemma}

\begin{proof}
We prove this standard statement on the corresponding tangent spaces.
\end{proof}

We will start with the case of the two-orbit $SL_2$-varieties.
We work out this case separately not only because it is simple
but also because it will be used quite often later on. 

\begin{proposition}
The two-orbit $SL_2$-varieties are $\mathbb P^1\times \mathbb P^1$
and $\mathbb P^2$ with the obvious $SL_2$-actions.
The respective generic stabilizers are the maximal torus 
$T$ and its normalizer.
\label{SL_2-var}
\end{proposition}

\begin{proof}
The a priori possible dimensions for a two-orbit $SL_2$-variety 
(that is, the possible dimensions for a $SL_2$-orbit)
are $0,1,2$ or $3$.

Using the fact that projective $SL_2$-orbits are 0- or 1-dimensional
(and vice versa) and that 3-dimensional $SL_2$-orbits are affine,
we prove that a two-orbit $SL_2$-variety can only be 2-dimensional.
Therefore, it has at least three points fixed by $T$ a maximal torus of
$SL_2$ (\emph{see}\/~\cite[Theorem 25.2]{Hum}) 
Since a projective $SL_2$-orbit has only two $T$-fixed
points, a generic stabilizer $H$ must contain this torus $T$.
Finally, we get that $H/T$ is finite because $\dim H=1=\dim T$ so
$H=T$ or $H=N_G(T)$.
The respective two-orbit $SL_2$-varities are clearly 
$\mathbb P^1\times \mathbb P^1$ and $\mathbb P^2$.
Note that the action of $SL_2$ on $\mathbb P^2$ is induced by $SO_3$.
\end{proof}

From now on, the group $G$ is of rank greater than $2$.
Let us start by studying locally the two-orbit $G$-variety $X$. 
Consider $y$ the $B$-fixed point of $X$ and $P$ its stabilizer in $G$. 
The group $P$ is parabolic and
$P^uL$ will be its Levi decomposition such that $T\subset L$. 
By~\cite[Theorem 1.4]{BLV}, there exists
an affine $L$-stable subvariety $Z$ of $X$ 
such that $\dim Z\geq 1$ and $Z\cap G\cdot
y=\{y\}$. 
Using the fact that $T$-orbits are affine,
it is easy to see that an affine $T$-variety containing fixed points
also contains $1$-dimensional orbits if the action is not trivial.
Therefore, there is a generic element $x$ of $X$ such that $\codim_T T_x\leq 1$.
In other words, we have just obtained:

\begin{proposition}
The ranks of the algebraic group $G$
and of the generic stabilizer differ at most by $1$.
\label{rang}
\end{proposition}

\begin{definition}
The two-orbit varieties such that the
rank of $G$ is equal to (resp. different from) the rank of $H$
are said of \emph{type~I} (resp. \emph{of type~II}).
\end{definition}

Suppose that the generic element $x$ is such that the torus
$T_x^\circ$ is a maximal torus of $G_x$. Let $L'$ be the
centralizer in $G$ of $T_x^\circ$. Then, by the above
proposition, $L'$ is equal to $T$ or is a Levi subgroup of
semisimple rank~$1$. We are going to prove that the latter
possibility can not occur \emph{i.e.}

\begin{proposition}
The maximal tori of the generic stabilizer are regular tori 
of the group $G$.
\label{reg}
\end{proposition}

\begin{proof}
We proceed by contradiction; suppose that $L'\neq T$ and
consider the $L'$-variety $\closure(L'\cdot x)$. 
Since $T_x^\circ\neq T$, 
$\closure(L'\cdot x)$ is actually a $SL_2$-variety but not a two-orbit $SL_2$-variety
(\emph{see}\/ Proposition~\ref{SL_2-var}). 
More precisely, 
$\closure(L'\cdot x)\setminus L'\cdot x$
is a finite union of closed $L'$-orbits of dimension $\leq 1$
all contained in $G\cdot y$ because of Lemma~\ref{compos}.
Therefore (\emph{see}\/ also the proof of Proposition~\ref{SL_2-var}),
the $L'$-variety $\closure(L'\cdot x)$ is of dimension $2$.
 
Let $B_{L'}$ be the Borel subgroup of $L'$ equal to $B\cap L'$ 
and $U_{L'}$ be its unipotent radical. 
There exists an element, say $x'$, in $L'\cdot x$ fixed by $U_{L'}$
(\emph{see}\/ the introduction). 
Moreover, $T_x^\circ$ fixes obviously $x'$ so, by Proposition~\ref{rang},
$T_x^\circ=T_{x'}^\circ$. 
Thus, 
$\closure(L'\cdot x')=
L\cdot\closure(B_{L'}\cdot x')=L'\cdot\closure(T\cdot x')$. 
If $y_1$ and $y_2$ denote
the $T$-fixed points in $\closure(T\cdot x')$ then 
$\closure(L'\cdot x)=L'\cdot x\cup L'\cdot y_1\cup L'\cdot y_2$.
The variety $\closure(L'\cdot x)$, being 2-dimensional, 
has at least three $T$-fixed points,
(\emph{see}\/~\cite[Theorem 25.2]{Hum}).
Therefore, $y_1$ and $y_2$ can not be simultaneously fixed by $L'$. 
Suppose that $L'\cdot y_1\neq y_1$ and consider
the element, say $y'$, satisfying $\closure(T\cdot y')=L'\cdot
y_1$.
Then we get:
$$
T_{y'}^\circ=T_{x'}^\circ, \quad
\dim(T\cdot x')=\dim(T\cdot y')=1 \quad \mbox{and} \qquad
y_1\in \closure(T\cdot x')\cap\closure(T\cdot y').
$$
It results
from the lemma below that $y_2\in L\cdot y_1$ which is
incompatible with the fact that $\closure(L'\cdot x)$ is not
a two-orbit variety. Thus the assumption $L'\neq T$ was
absurd -- which proves Proposition~\ref{reg}.
\end{proof}

\begin{lemma}
Let $x$ and $z$ be two elements of a
projective $G$-variety such that
$$
\dim T\cdot x=\dim
T\cdot z=1 \qquad \mbox{and} \qquad T_x^\circ=T_z^\circ.
$$
There are two or four $T$-fixed points in 
$\closure(T\cdot x)\cup \closure(T\cdot z)$ 
assuming that they are all contained in one same closed $G$-orbit.
\label{2or4}
\end{lemma}

\begin{proof}
By embedding the given variety $G$-equivariantly into the projective
space of a representation, we can write the elements $x$ and $z$ as:
$x=[v_{\mu_1}+\ldots+v_{\mu_r}]$ and 
$z=[v_{\nu_1}+\ldots+v_{\nu_s}]$.
By assumption, the weights of the support of $x$ 
(resp. of the support of $z$)
sit on a same affine line $D_x$ (resp. $D_z$) and moreover,
$$
\ker(\mu_1-\mu_r)^\circ=T_x^\circ=T_z^\circ=\ker(\nu_1-\nu_s)^\circ.
$$
From these equalities, 
we can deduce that the lines $D_x$ and $D_z$ are strictly parallel or equal.
The first situation clearly yields four $T$-fixed points in 
$\closure(T\cdot x)\cup \closure(T\cdot z)$
and the second one only two $T$-fixed points (use the fact that 
the Weyl group acts transitively on 
the $T$-fixed points of a projective orbit).
\end{proof}

As a consequence of Proposition~\ref{reg}, we have:

\begin{corollary}
Consider two elements $x$ and $z$ of the
dense $G$-orbit of a two-orbit variety such that
$T_x^\circ$ (resp. $T_z^\circ$) is a maximal torus of $G_x$
(resp. $G_z$). Then $z\in N_G(T)\cdot x$; in particular, if
$x$ and $z$ are $B$-conjugated they are actually $T$-conjugated.
\label{conjug}
\end{corollary}

\section{Parabolic induction}
\label{induction}

This section is the main step to the classification of
two-orbit varieties. It consists essentially in proving
that the problem of classifying the two-orbit varieties
can be reduced to a subclass of two-orbit varieties called
\emph{cuspidal}, \emph{i.e.} the two-orbit varieties which
can not be obtained by parabolic induction 
(\emph{see} section~\ref{definduc}).

\subsection{Statements}

\begin{proposition}
If $X$ is obtained by parabolic induction from a pair $(P,Y)$,
then it is completely determined by its closed orbit and the
variety $Y$.
\label{geometry}
\end{proposition}

\begin{proof}
This statement is a direct consequence of the
fact that the morphism $G\times_P Y\rightarrow X$ is a
blowing down of the closed $G$-orbit of $X$ and that the
variety $X$ is completely determined, according to a result
of \textsc{Luna} and {Vust} (\emph{see}\/~theorem~8.3
of~\cite{LV}), by its $G$-stable prime divisors and by the
$B^-$-stable prime divisors of its dense $G$-orbit, whose
closure in $X$ contains the closed $G$-orbit of $X$.
\end{proof}

We can put an order $\leq$ on the set of induction pairs of
$X$, defined naturally, for two pairs
$(P_1,Y_1)$ and $(P_2,Y_2)$, by:
\begin{equation*}
(P_1,Y_1)\leq (P_2,Y_2) \qquad \mbox{if $P_1\subseteq P_2$
and $\varphi_1(Y_1)\subseteq \varphi_2(Y_2)$}.
\end{equation*}

\begin{theorem}
The set of induction pairs of a two-orbit
variety (endowed with the order~$\leq$) has an unique
minimal element.
\label{induc}
\end{theorem}

\subsection{Proof of Theorem~\ref{induc}}

\label{proof}

In order to prove Theorem~\ref{induc}, we are
constructing explicitly, in this section, the minimal
induction pair.

Let us consider an induction pair $(P,Y)$ and two elements
$x$ and $z$ in the dense $G$-orbit of $X$ such that $B\cdot
x$ and $B\cdot z$ are closed in $G\cdot x$. These elements
can be chosen such that $T_x^\circ$ (resp. $T_z^\circ$) is
a maximal torus of $G_x$ (resp. $G_z$).
Then, by Corollary~\ref{conjug}, 
they are conjugated by an element $n\in N_G(T)$. 
Because of the choice made on $B\cdot x$
and on $B\cdot z$, one can show by a short computation that $n$ is 
in fact an element of $P$.
We have got:

\begin{lemma}
If $X$ is obtained from a pair $(P,Y)$ by
parabolic induction then the $B$-orbits which are closed in
the dense $G$-orbit are in a same $P$-orbit of $X$. In
particular, they are contained in $\varphi(Y)$.
\end{lemma}

Let $x$ be as above, denote $G_x$ by $H$. 
Let $P_1$ be the parabolic subgroup of
$G$ generated by $B$, $H$ and the elements $n$ of $N_G(T)$
such that $Bn\cdot x$ is closed in $G\cdot x$. By
construction, for all pairs $(P,Y)$, we have: $P_1\subset
P$ and $X_1=\closure(P_1\cdot x)\subset \varphi(Y)$. 
However, $X_1$ may not be a two-orbit $P_1$-variety. 
More precisely, we have:
$X_1\cap G\cdot x=P_1\cdot x$ 
but an analogous equality may not hold for
$X_1\cap G\cdot y$ ($y$ a $T$-fixed element of the closed
orbit in $X_1$). 
So, instead of $P_1$, we have to consider
$P_H$ the parabolic subgroup of $G$ generated by $P_1$ and
the elements $w\in W$ such that $w\cdot y\in X_1$. Thus,
$X_H=\closure(P_H\cdot x)=\closure(P_H\cdot X_1)$ is a two-orbit
$P_H$-variety with $P_H\subset P$ and $X_H\subset
\varphi(Y)$.

Finally, let us consider $\widetilde X_H$ the normalized
variety of $X_H$.
The pair $(P_H,\widetilde X_H)$ is the required element.

\begin{proposition}
The pair $(P_H,\widetilde X_H)$ is an
induction pair; this is the minimal element for the set of
induction pairs of $X$.
\label{minimal}
\end{proposition}

We already know that the variety $\widetilde X_H$ is a
two-orbit $P_H$-variety.
So we are left only to prove that
the radical of $P_H$ acts trivially on $\widetilde X_H$.
For this, we use the following

\begin{proposition}
Let $Y$ be a two-orbit $Q$-variety
with $Q$ a parabolic subgroup of $G$ containing $B$.
Consider an element $z$ of its dense $Q$-orbit such that
the torus $T_z^\circ$ is a maximal torus of $Q_z$. 
Suppose that the radical of $Q$ doesn't act trivially on $Y$.
Then the orbit $L\cdot z$ is complete for $L$ a Levi subgroup of $Q$
containing the torus $T$.
\label{radical}
\end{proposition}

\def\prooftext{of Proposition~\ref{minimal}}
\begin{proof}
According to Proposition~\ref{radical} applied to
$\widetilde X_H$ and $P_H$, 
if we prove that $L\cdot z$ is not complete,
Propostion~\ref{minimal} will follow. 
If $T_z^\circ\neq T$, that $L\cdot z$ is not complete is obvious. 
Consider then the other case: $T_z^\circ=T$. 
If $\alpha$ is a simple root of $P_H$ such
that $U_\alpha\not\subset G_z$ then (because of the
assumption made on $B\cdot z$) $s_\alpha Bs_\alpha\cap G_z$
contains a Borel subgroup of $G_z$. In other words,
$s_\alpha Bs_\alpha$ is closed in $P_H\cdot z$ and by
definition of $P_H$, it means that $SL_2(\alpha)\subset
L$. To conclude that $L\cdot z$ is not complete, 
we need only to observe that $\closure(L\cdot z)$ contains a $T$-fixed
point of the closed $G$-orbit 
(because $\closure(U_\alpha\cdot z)$ does by Corollary~\ref{conjug}).
\end{proof}

\def\prooftext{of Proposition~\ref{radical}}
\begin{proof}
First of all, note that the radical of $Q$ acts trivially
on the closed $Q$-orbit of $Y$.
So the proposition relies only on the dense $Q$-orbit.

We choose an element $z$ verifying the hypotheses of the
proposition and such that $B\cdot z$ is closed in
$Q\cdot z$.

Let us start with the case: $T_z^\circ\neq T$.
Thus $L\cdot z$ is not complete. 
Let $\beta\in \Phi^+$ be such that $U_\beta\not\subset Q_z$. 
Such a root exists otherwise $Y$ will have two $B$-fixed points 
given by the $T$-fixed points of $\closure(T\cdot z)$. 
Consider the $2$-dimensional variety $X_\beta=\closure(U_\beta T\cdot z)$ and denote by
$T_{X_\beta}$ ($\subset\ker\beta$) its generic stabilizer
in $T$. To get our result, we are going to prove that
$T_{X_\beta}^\circ$ contains the identity-component
$Z(L)^\circ$ of the center of $L$.

Let $y_1$ and $y_2$ be the $T$-fixed points of 
$\closure(T\cdot z)$. 
If $y_i$ ($i=1,2$) is not fixed by $U_\beta$, 
we denote by $y_i'$ the other $T$-fixed point in
$\closure(U_\beta\cdot y_i$); otherwise, we set $y_i'=y_i$. 
The points $y_1'$ and $y_2'$ are distinct; 
this follows easily from Lemma~\ref{2or4} and 
from the fact that $T_z^\circ$ is regular.

The variety $X_\beta^{U_\beta}$, being connected and
containing two distinct elements, $y_1$ and $y_2$, is $1$-dimensional.

\begin{claim}
There exist $z_1,\ldots,z_r$ and
$\rho_1,\ldots,\rho_{r+1}$ ($r\geq 1$) in
$X_\beta^{U_\beta}$ such that $\dim T\cdot z_i=1$ and such that
$\rho_i,\rho_{i+1}$ are the $T$-fixed points in $\closure
(T\cdot z_i)$. Moreover, if $y_i=y_i'$ ($i=1,2$) then
$r\geq 2$.
\end{claim}

To construct the elements $z_i$ and $\rho_i$, it suffices
to consider the convex hull of the support of $X_\beta$;
the elements $z_i$ (resp. $\rho_i$) have as support the
edges (resp. vertices) of this polytope.

Becasuse of the claim, we have now at hand 
at least two elements $u$ and $v$
simultaneously in the closed $Q$-orbit and in $X_\beta$
such that $\dim(T\cdot u)=1=\dim(T\cdot v)$ and
$\closure(T\cdot u)\cap \closure(T\cdot v)=1$.
Therefore:
$\bigl(T_u\cap T_v\bigr)^\circ$ (of codimension $2$) is
equal to $T_{X_\beta}^\circ$.
But since $Z(L)$ acts
trivially on the closed $Q$-orbit, we end up with the required inclusion:
$Z_L\subset T_{X_\beta}^\circ$.

\hangindent=0pt \hangafter=1

Assume now that $T_z^\circ=T$ and that $L\cdot z$ is not
closed. Therefore, there exists at least one positive root,
say $\alpha$, in the root system of $(L,L)$ such that
$U_\alpha\not\subset Q_z$. 
Recall that we want to prove that the radical of $Q$ acts trivially on
$Y$.
To do so, we proceed by contradiction:
suppose there exists
$\beta\in\Phi^+$, $\beta\neq\alpha$ such that
$U_\beta\subset Q^u$ and $U_\beta\not\subset Q_z$. 
Then, the variety $X_{\beta,\alpha}=\closure(U_\beta U_\alpha\cdot
z)$ is 2-dimensional and contains a dense $T$-orbit; denote
by $T_{X_\beta}$ its generic stabilizer in $T$.


If $y_1$ and $y_2$ are the $T$-fixed points (distinct from
$z$) of $\closure(U_\alpha\cdot z)$ and of $\closure(U_\beta\cdot z)$ respectively, 
then necessarily $y_1$ and $y_2$ are distinct (as they do not have
the same support). 
If $U_\beta\cdot y_1\neq y_1$, we denote
by $y_1'$ the other $T$-fixed point of $\closure(U_\beta\cdot
y_1)$; otherwise, we set $y_1'=y_1$. 
Then the  variety $X_{\beta,\alpha}^{U_\beta}$ is 1-dimensional,
since it contains the distinct points $y_1'$ and $y_2$. 
But, $y_1$ and $y_2$ must be the only $T$-fixed points of the closed
$Q$-orbit of $Y$ in $X_{\beta,\alpha}$ otherwise with the
same arguments used in the first case ($T_z^\circ\neq T$), 
we will have: $Z(L)^\circ\subset T_{X_\beta}$. 
So, we can conclude that $U_\beta$ must fix $y_1$ 
and that there exists an element $u$ in the closed $Q$-orbit of $Y$ such that:
$X_{\beta,\alpha}^{U_\beta}=\closure(T\cdot u)$.

By considering the variety $\closure(U_{-\alpha}\cdot y_1)$, we
get $y_3$, the other $T$-fixed point in it.
For the same reasons as before, the points $y_1$ and $y_3$ are distinct.

To conclude, we have constructed two 1-dimensional subvarieties of the
closed $Q$-orbit, $X_{\beta,\alpha}$ and
$\closure(U_{-\alpha}\cdot y_1)$ such that $y_1$ is their only common $T$-fixed point. 
According to Lemma~\ref{2or4}, this yields the contradiction:
$(T_u\cap\ker\alpha)^\circ=T_{X_\beta}^\circ$.
\end{proof}

\section{Two-orbit varieties of type~I}
\label{typeI}

In this section, the two-orbit $G$-variety 
$X=\closure(G\cdot x)=G\cdot x\cup G\cdot y$
is of type I, \emph{i.e.} $T_x^\circ =T_x$.
Recall (\emph{see}~section~\ref{notation}) that a two-orbit variety is
projective.
So we can embed $X$ in $\mathbb P(V)$ with $V$ a finite $G$-module.
The elements $x$ and $y$ of $X$ can be written as:
$x=[v_\mu]$ and $y=[v_\lambda]$
with $v_\mu$ and $v_\lambda$ weight vectors of $V$.
We can choose $\lambda$ to be dominant.

Let $\beta\in\Phi^+$ be such that $U_\beta\not\subset G_x$
(such a root exists otherwise $B\subset G_x$) 
and consider
$z=[v_\mu+v_{\mu+\beta}\cdots+v_{\mu+k\beta}]\in U_\beta\cdot x$. 
With a  judicious choice of $x$, (for instance $\mu$ dominant),
one shows easily that $[v_{\mu+k\beta}]\in G\cdot y$.
So finally, we can choose $x$ such that $\lambda=\mu+k\beta$
(take a $W$-conjugate of the previous $x$).
Define the support of a root (denoted $\supp $) as the set of simple
roots really involved in its writting.

\begin{proposition}
If $X$ is cuspidal then $G$ is simple
and $\supp \beta=\Delta$.
\label{cuspid1}
\end{proposition}

\begin{proof}
The first assertion comes from the
construction of the minimal induction pair 
(\emph{see}~\/section~\ref{proof}). 
More precisely, as $X$ is cuspidal,
$G$ must be generated by the parabolic subgroup $P_1$ and
the elements $w\in W$ such that $w\cdot[v_\lambda]\in X_1$.
Recall that $P_1$ is spanned by $H$, $B$ and the elements
$n\in N_G(T)$ such that $Bn\cdot x$ is closed in $G\cdot x$
and $X_1=\closure(P_1\cdot x)$.

If $G=G_1\times\cdots\times G_r$ with $G_i$ simple and
$\Delta=\Delta_1\times\cdots\times \Delta_r$ with
$\Delta_i$ associated to $G_i$, we are going to show that
$H=H_1\times\cdots\times H_r$ with $H_i=G_i$ for all $i\neq
i_0$ and $\supp\beta\subset\Delta_{i_0}$.

First of all, note that 
if $\alpha\not\in\supp\beta$, then $U_\alpha\subset G_x$.
Indeed, $\mu+\ell\alpha$ ($\ell >0$) is not of shape
$\lambda-\sum_{\gamma\in\Delta}n_\gamma\cdot\gamma$,
(weights of $V$) since $\lambda=\mu+k\beta$ 
(\emph{see} the choice of $x$ made above).
Moreover, if $\alpha\in\Delta\setminus\Delta_{i_0}$ then
with the above description of $G$, we must have
$U_{-\alpha}\subset G_x$. 
Thus: $U_{\pm\alpha}\subset G_x$,
for all $\alpha\in\Delta\setminus\Delta_{i_0}$.
The acting group $G$ can then be assumed to be simple.

To obtain the second assertion, consider the parabolic
subgroup $P$ associated to $\supp\beta$ and the $P$-variety
$Z=X \cap \mathbb P \bigl(\bigoplus_\nu V_\nu \bigr)$ for
$\nu=\lambda-\sum n_\alpha\alpha$ with $\alpha \in
\supp\beta$ and $n_\alpha\geq 0$. Then, if
$\supp\beta\neq\Delta$, $(P,Z)$ is an induction pair of
$X$.
\end{proof}

From now on, $G$ will be simple and $x$ will satisfy the
above conditions as well as the two following ones. If
$\lambda=\mu+k\beta$ then $k$ is minimal and moreover if
$[v_{\mu'}]$ is another $T$-fixed generic element
satisfying the same assumptions as $x$ then $\mu>\mu'$ for
$\mu$ and $\mu'$ comparable.

Consider a second positive root, say $\gamma$, with
$\gamma\neq\beta$ and $U_\gamma\not\subset G_x$ 
(there exists at least one such a root which is simple).
Denote by $L$ the Levi subgroup
associated to $\beta$ and $\gamma$, \emph{i.e.} 
$L$ is the centralizer in $G$ of $(\ker\beta\cap\ker\gamma)^\circ$.
Then we have:

\begin{lemma}
The variety $\closure(L\cdot x)$ is a cuspidal two-orbit
$L$-variety.
\end{lemma}

\begin{proof}
From Lemma~\ref{compos}, we have $L\cdot
x=\closure(L\cdot x)\cap G\cdot x$ and 
$\closure(L\cdot x)\setminus L\cdot x$ 
consists of a finite union of complete $L$-orbits. 
Thus, if $\closure(L\cdot x)$ is not a two-orbit $L$-variety, 
we will have $(L\cdot x)^{U_L}\neq\emptyset$ for
$U_L$ an unipotent maximal subgroup of $L$ 
(again the same argument as in the introduction).
It will follow that $\closure(L\cdot x)=L\cdot x$. 
But this equality can not hold
since $[v_\lambda]\in\closure(L\cdot x)\setminus L\cdot x$.

As we can not find any proper parabolic subgroup $P$ such
that $P^u\subset L_x\subset P$ (because $U_{\pm\beta}$,
$U_{\pm\gamma}\not\subset G_x$), $L\cdot x$ must be
cuspidal.
\end{proof}

By this procedure, we have constructed two-orbit varieties
for some subgroups $L$ of $G$ of semisimple rank~$2$. Thus, once
we know what the two-orbit varieties are, for the simple
groups of rank~$2$, we will know $L_x$ (for all $L$'s) and
then $G_x$.

Let us start with determining the cuspidal two-orbit
varieties in the rank~$2$ case. For this, we need two
technical lemmas.
Let $\Lambda$ be the convex hull of the support of the
variety $\closure(U_\gamma U_\beta\cdot x)$. 
We have the following picture and notation in the weight lattice
$\mathcal X$ of $G$.

\begin{lemma}
The extremal points $\nu_i$ of $\Lambda$ are
$W$-conjugated to $\lambda$ and the $\beta_i$'s are some
roots of~$\Phi$.
\label{hull}
\end{lemma}

\begin{proof}
First of all, note that the points $[v_{\nu_i}]$ 
and the elements $[v_i]$ whose support is
$[\nu_i,\nu_{i+1}]\cap \mathcal X$ 
sit in $\closure(U_\gamma U_\beta\cdot x)$. 
If $\Lambda$ has more than three extremal
points then the cardinality of $\supp [v_i]$ is smaller
than $\supp x$'s one.
So by the minimality assumption made on $x$, we must have
$[v_i]\in G\cdot y$ and thus $[v_{\nu_i}]\in G\cdot y$. 
If $\Lambda$ has three extremal points, 
the assumed maximality of the weight $\mu$ ($x=[v_\mu]$) forces 
$\nu_1$ to be in $W\cdot \lambda$.

The second assertion of the lemma follows from the fact
that $T_{[v_i]}^\circ=\ker(\nu_i-\nu_{i+1})^\circ$ is singular.
\end{proof}

\begin{lemma}
Let $\delta$ be a positive root such that
$[v_{s_{\gamma'}(\lambda)}]\in \closure(U_\delta\cdot x)$, for
$\gamma'\in\Delta$. Suppose that: if there exists $r\geq 0$
such that $\mu+r\gamma'$ is extremal as weight of $V$, we
must have $\mu+r\gamma'=s_\alpha(\lambda)$ for
$\alpha\in\Delta$. Then, $U_{\gamma'}\subset G_x$ if
$(\mu,\gamma')\geq 0$ and $U_{-\gamma'}\subset G_x$ if
$(\mu,\gamma')\leq 0$.
\label{critere}
\end{lemma}

\begin{proof}
If there is no $r>0$ such that $\mu+r\gamma'$ is
extremal then by Lemma~\ref{hull}, we must have
$U_{\gamma'}\subset G_x$ if $(\mu,\gamma')\geq 0$. If
$(\mu,\gamma')\leq 0$ and $U_{-\gamma'}\not\subset G_x$
then $\closure(U_{-\gamma'}\cdot x)$ contains a $T$-fixed point
of $G\cdot y$ and so does $\closure(U_{\gamma'}\cdot x)$ -- which is absurd.

Suppose now that $\mu+r\gamma'$ is extremal and that
$(\mu,\gamma')\geq 0$. 
Let $x'=\varepsilon_\delta(1)\cdot x\in U_\delta\cdot x$,
$\varepsilon_\delta$ being the natural map associated to $\delta$ 
from $\mathbb C$ to $U_\delta$. 
Then by assumption
$x'=[v_\mu+\cdots+v_{s_{\gamma'}(\lambda)}]$. 
The study of
$\supp(\varepsilon_\beta(1)\cdot x')$ gives raise to a $j>0$ such that
\begin{equation}
-\frac{Y_{\beta}^{n_{j}}}{n_{j}!}\cdot
v_{\mu +j\delta}=sY_{\gamma'}\cdot
v_{s_{\gamma'}(\lambda)}
\qquad \mbox{for $s\in \mathbb C\setminus\{0\}$ and
$n_j\geq 0$.}
\label{gap}
\end{equation}

\noindent
Recall (\emph{see}~\/section~\ref{notation}) 
that $Y_\beta$ denotes an element of $\mathfrak g_\beta\setminus\{0\}$.

Furthermore, consider the support of the variety $
\closure\bigl(T\cdot\exp (Y_{\beta}+sY_{\gamma'})\cdot x'\bigr)$.
The points $s_\alpha(\lambda)$, $s_{\gamma'}(\lambda)$ are extremal
points of this support but $\lambda$ is not; 
otherwise, we will have a gap in this support since the weight
$\gamma'+s_{\gamma'}(\lambda)$ is missing by
equality~(\ref{gap}). It implies that
$[s_\alpha(\lambda),s_{\gamma'}(\lambda)]$ must be an edge
of this support -- which contradicts Lemma~\ref{hull}
since $s_\alpha(\lambda)-s_{\gamma'}(\lambda)$ is not a
root, $\alpha$ and $\gamma'$ being simple.

If $(\mu,\gamma')\leq 0$, we can go back to the positive
case with the element $s_{\gamma'}\cdot x$.
\end{proof}

Let us show how we can apply these two lemmas to get the
two-orbit varieties in the rank~$2$ case, through the
following

\begin{example}
Suppose $G$ of type $B_2$. Then there are
two possibilities for $\beta$: $\beta=\alpha_1+\alpha_2$ or
$\beta=\alpha_1+2\alpha_2$. 
Let us compute $\mathfrak g_x$ in case $\beta=\alpha_1+\alpha_2$.

%

If $U_{\alpha_1}\not\subset G_x$ 
then $\alpha_2$ verifies the conditions of Lemma~\ref{critere} 
with $(\mu,\alpha_2^\vee)\geq 0$.
Therefore $U_{\alpha_2}\subset G_x$. Let us show that
$U_{\alpha_1+2\alpha_2}\subset G_x$. If there does not
exist $r>0$ such that $\mu+r(\alpha_1+2\alpha_2)$ is extremal
then this inclusion is given by Lemma~\ref{hull}.
Otherwise, we will have: $\mu+r\alpha=s_1(\lambda)$. A simple
computation leads to:
$2(\lambda,\alpha_1^\vee)=(\lambda,\alpha_2^\vee)$. 
But in this latter case, $\mu$ does not satisfy the good conditions;
in particular, ``$k$'' is not minimal. 
We have obtained the pair: 
$(SO_5,\mathfrak g_x=\mathfrak t\oplus\mathfrak
g_{\alpha_2}\oplus\mathfrak g_{\alpha_1+2\alpha_2})$.


If $U_{\alpha_1}\subset G_x$ then $U_{\alpha_2}\not\subset
G_x$ and we have $(\mu,\pm\alpha_1)=0$.
Thus by Lemma~\ref{critere}, we get
$U_{\pm\alpha_1}\subset G_x$. 
It follows from Lemma~\ref{hull} that
$U_{\alpha_1+2\alpha_2}\subset G_x$.
Therefore, we have obtained the pairs:
$(SO_5,\mathfrak{so}_4)$ and $(SO_5,\mathfrak g_x=\mathfrak
t\oplus\mathfrak g_{\pm\alpha_1}
                        \oplus\mathfrak
                        g_{\alpha_1+2\alpha_2})$.
\end{example}

By the procedure given in this example, we get the
two-orbit varieties for the rank~$2$ case.
In other words, we have a description of $\Phi(G_x)$ the root system of
$G_x$.

\begin{lemma}
Let $\gamma\in\Phi^+$, we have the following:

\smallbreak{\rm(i)}\enspace
if $\langle \gamma
,\beta \rangle $ is of type $A_1\times A_1$ then
$\pm\gamma\in\Phi (G_x)$;

\smallbreak{\rm(ii)}\enspace
if $\langle \gamma,\beta \rangle
$ is of type $A_2$ then $\pm \gamma \in \Phi (G_x)$ or $\pm
s_\gamma(\beta) \in \Phi (G_x)$;

\smallbreak{\rm(iii)}\enspace
if $\langle \gamma,\beta \rangle
$ is of type $B_2$ with $\beta =\varepsilon _1+\varepsilon
_j$ then $\Phi(G_x)\cap \langle \gamma,\beta \rangle
=\{\varepsilon_1-\varepsilon _j,\varepsilon_1\}$;

\smallbreak{\rm(iv)}\enspace
if $\langle \gamma,\beta \rangle
$ is of type $B_2$ with $\beta =\varepsilon_1$ then $\Phi
(G_x)\cap\langle \gamma,\beta \rangle
\supset\{\varepsilon_1+\varepsilon_j,\pm
(\varepsilon_1-\varepsilon _j)\}$ or $\Phi (G_x)\cap\langle
\gamma,\beta \rangle =
\{\varepsilon_1+\varepsilon_j,\varepsilon_j\}$.
\label{reformul}
\end{lemma}

By applying this statement to $\beta$, a positive root of maximal
support and to any other positive root, 
we get Table~1 just by computations. 
Let us work out two examples to understand how it works.

\begin{example}
Suppose $G$ of type $A_n$. 
Then $\beta=\alpha_1+\cdots+\alpha_n$. 
Since
$\langle\alpha_i,\beta\rangle$ is of type $A_1\times A_1$,
by Lemma~\ref{reformul}-(i), 
we have $\pm\alpha_i\in\Phi(G_x)$ for $2\leq i\leq n-1$. 
Moreover, Lemma~\ref{reformul}-(ii) yields:
$\pm\alpha_1\in\Phi(G_x)$ or $\pm\alpha_n\in\Phi(G_x)$. 
This gives the pair $(SL_{n+1},\mathfrak{gl}_n)$.

Suppose $G$ of type $B_n$ and $\beta=\varepsilon_1$. 
Applying Lemma~\ref{reformul}-(i) and Lemma~\ref{reformul}-(ii) respectively, 
we get
$\{\pm\alpha_2,\ldots,\pm\alpha_n\}\subset\Phi(G_x)$ 
and respectively
$\{\varepsilon_{1}+\varepsilon_{i}:i\geq 2\} \subset \Phi (G_{x})$ 
with
$\pm(\varepsilon_1-\varepsilon_i)\in\Phi(G_x)$ or
$\varepsilon_i\in\Phi(G_x)$ for all $i\geq 2$. 
If $\pm\alpha_1\not\in\Phi(G_x)$ then 
$\{\pm\alpha_1,\ldots,\pm\alpha_{n-1},\alpha_{n-1}+2\alpha_n\}\subset\Phi
(G_x) $. 
Therefore, we obtain the two pairs
$(SO_{2n+1},\mathfrak{so}_{2n})$ and 
$\bigl(SO_{2n+1},\mathfrak{g}_x=\bigoplus_{\alpha\in\Psi}\mathfrak
g_\alpha \oplus\mathfrak t \bigr) $ where
$\Psi=\langle\pm\alpha_1,\ldots,\pm\alpha_{n-1},\alpha_{n-1}+2\alpha_n\rangle$.
If $\pm\alpha_1\not\in\Phi(G_x)$, we get the pair
$(SO_{2n+1},\mathfrak{g}_x
        =\bigoplus_{\alpha\in\Psi}\mathfrak
        g_\alpha\oplus\mathfrak t)$, with $\Psi=\langle \pm
\alpha_2,...,\pm
\alpha_{n-1},\alpha_n,\varepsilon_1+\varepsilon_n\rangle$.
\end{example}

\section{Two-orbit varieties of type~II}
\label{typeII}

In this section, the two-orbit $G$-variety 
$X=\closure(G\cdot x)$ is of type~II, that is, by definition,
$T_x^\circ\neq T$.
We embed $X$ in $\mathbb P(V)$ with $V$ a finite $G$-module as in the
previous section.
Then the generic element $x$ can be written as
$[v_{\lambda_0}+\cdots+v_{\lambda_{r+1}}]$ with the
$\lambda_i$'s sitting on a same affine line, say $D_x$
(because of Proposition~\ref{rang}).
we order the weights $\lambda_i$ of $\supp x$ in such way that:
$\closure(T\cdot x)= T\cdot x\cup
\left\{[v_{\lambda_0}],[v_{\lambda_{r+1}}]\right\}$.
As elements of $G\cdot y$,
$[v_{\lambda_0}]$ and $[v_{\lambda_{r+1}}]$ are $W$-conjugate, \emph{i.e.}
there exists $w\in W$ such that $\lambda_0=w(\lambda_{r+1})$.
We choose $x$ such that $\lambda_0=\lambda$
and satisfying the following condition of minimality:
if $x'=[v_\lambda+\cdots+v_{w'(\lambda)}]$ is another generic element then
$w<w'$ for $w$ and $w'$ comparable.

Let $\alpha$ be a simple root such that $s_\alpha w<w$ and
$U_\alpha\not\subset G_x$.
Take for instance $\alpha$ such that $\bigl(w(\lambda),\alpha\bigr)<0$. Consider the
variety $X_\alpha = \closure(TU_\alpha \cdot x)$ and in
particular, the convex hull $\Lambda$ of its support
pictured below in the weight lattice $\mathcal X$ of $G$.


Similarly as for Lemma~\ref{hull}, we get:

\begin{lemma}
The elements $y_i\in X_\alpha$ whose support
sits on the line $D_i$ belong to the closed $G$-orbit of
$X$. The directions of the affine lines $D_i$ are given by
roots $\beta_i$. In particular, the extremal points
$\nu_i$ are $W$-conjugated and $\supp
y_i=[\nu_i,\nu_{i+1}]\cap \mathcal X$.
\end{lemma}


The roots $\beta_i$ span a root system of rank~$2$; let
$\{\alpha,\beta\}$ be a basis of this root system.

\begin{corollary}
$w=s_\alpha s_\beta$.
\end{corollary}

\begin{proof}
We know that $s_\alpha w<w$ and that $\lambda
-w(\lambda)$ can not be, up to a scalar, a root because
$T_x^\circ=\ker\bigl(\lambda -w(\lambda)\bigr)^\circ$ is
regular (\emph{see}~Proposition~\ref{reg}). 
So we can assume $\langle\alpha,\beta\rangle$ to
be of type $G_2$ with $w=(s_\alpha s_\beta)^2$~-- the case
$w=w_0$ being easily ruled out.

Consider the convex hull $\Lambda$ 
(\emph{see}~\/ the corresponding picture). 
Since we have: $\supp y_i=[\nu_i,\nu_{i+1}]\cap \mathcal X$,
there exists a weight $\nu\in D_x$ such that $\lambda
-\beta_t=\nu+k\alpha$, $k\geq 0$. In other words, if
$D_{\alpha,\lambda-\beta_t}$ denotes the line of direction
$\alpha$ passing through $\lambda-\beta_t$, we must have:
\begin{equation}
D_x\cap D_{\alpha,\lambda-\beta_t}\cap
\mathcal{X}\neq \emptyset.
\label{intersect}
\end{equation}
If $w=(s_\alpha s_\beta)^2$ with
$\langle\alpha,\beta\rangle$ of type $G_2$ then
$\beta_t=\beta$ or $\beta=\beta+\alpha$. 
But (\ref{intersect}) forces
$(\lambda,\alpha)/\bigl((\lambda,\alpha)+(\lambda,\beta)\bigr)$
to be an integer -- which can not occur because
$(\lambda,\alpha)\cdot(\lambda,\beta)\neq 0$ since
$T_x^\circ$ is regular.
\end{proof}


Using this corollary and the same arguments given in the
proof of Proposition~\ref{cuspid1}, we get:

\begin{proposition}
If $X$ is cuspidal of type II then $G$
is simple or of type $A_1\times A_1$. Furthermore,
$\supp\beta\cup\{\alpha\}=\Delta$.
\label{rk2}
\end{proposition}

From now on $G$ will be assumed to be simple or of type
$A_1\times A_1$. 
Let $L$ be the Levi subgroup associated to $\alpha$ and $\beta$
and $\mathfrak l$ be its Lie algebra.

\begin{proposition}
{\rm(i)}\enspace
If $\langle\alpha,\beta\rangle$ is of type $A_1\times A_1$
then $\mathfrak l_x= \mathfrak t'\oplus \mathbb
C(Y_{-\alpha}+Y_\beta) \oplus \mathbb
C(Y_\alpha+Y_{-\beta})$.
\smallbreak
{\rm(ii)}\enspace
Otherwise, $\mathfrak l_x=
\mathfrak t'\oplus \mathbb C(Y_{-\alpha }+Y_{s_\alpha
(\beta )}) \oplus \bigoplus_{\substack{ \gamma \in \Phi
^+\cap\langle \alpha,\beta\rangle \\ \gamma \neq
\alpha,s_\alpha (\beta )}}\mathfrak{g}_\gamma $.
\smallbreak
\noindent
Here, $\mathfrak t'$ is the kernel of
$\alpha +s_\alpha (\beta)$ considered as element of the
dual $\mathfrak t^*$.
\label{listII}
\end{proposition}

\begin{proof}
Set $\delta=s_\alpha(\beta)$, $(s_\alpha
s_\beta(\lambda),\delta^\vee)<0$.
Thus $Y_\delta\notin
\mathfrak g_{[v_{w(\lambda)}]}$ and $Y_\delta\not\in
\mathfrak g_x$.

In order to get: $Y_\delta +Y_{-\alpha}\in \mathfrak g_x$,
we are going to prove
\begin{equation}
\label{egal}
(\lambda,\alpha^\vee)=(\lambda,\beta^\vee).
\end{equation}

Set $\lambda=m\omega_\alpha+n\omega_\beta$ ($m,n>0$) and
consider the variety $\closure(TU_\alpha\cdot x)$. 
The arguments used in the proof of the previous corollary give
a weight $\lambda_i\in\supp x$ such that $Y^r_\alpha \cdot
v_{\lambda_i}$ is of weight $\lambda-\beta$ (for $r>0$).
Translating the latter in terms of equations, we get
$\lambda_i=\lambda-\alpha_n/m\delta$ and $n/m\in \mathbf N$. 
Similarly, considering the variety $\closure(U_\alpha\cdot
x)$, we get $m/n\in \mathbf N$ thus $m=n$ and also
\begin{equation}
Y_{-\alpha}\cdot v_\lambda=-q
Y_\delta\cdot v_{\lambda_i} \quad q\in\mathbf C^*.
\label{gap2}
\end{equation}

Finally, let $Z_t$ be the variety $\closure\bigl(T\exp
t(Y_{\delta }+qY_{-\alpha })\cdot x\bigr)$, $t \in \mathbf
C$. Its support is entirely contained in the triangle of
vertices $\lambda $, $s_{\alpha }(\lambda )$ and $w(\lambda
)$. But there is a gap in this support: the weight
$\lambda-\alpha$ is missing because of (\ref{gap2}). This
implies that $\exp t(Y_{\delta }+qY_{-\alpha })\in G_x$.

To show that $U_\gamma\subset G_x$ for all
$\gamma\in\Phi\setminus\{\alpha,\delta\}$, we proceed by
contradiction and as before, we will find a gap in the
support of $\closure(U_\gamma\cdot x)$.

To conclude, we have to note that if
$\langle\alpha,\beta\rangle$ is not of type $A_1\times A_1$,
then $Y_\alpha+Y_{-s_\alpha(\beta)}\not\in\mathfrak g_x$.
The first assertion is obtained just by symmetry.
\end{proof}

As a consequence of the previous proof, we have
\begin{corollary}
{\rm(i)}\enspace
$(\lambda,\alpha^\vee)= (\lambda,\beta^\vee)$;
\smallbreak
{\rm(ii)}\enspace
$T_x^\circ =\ker \bigl(\alpha +s_\alpha (\beta
)\bigr)^\circ$;
\smallbreak
{\rm(iii)}\enspace
$N_G(H)/H$ is finite;
\smallbreak
{\rm(iv)}\enspace
$\closure(L\cdot x)$ is a cuspidal two-orbit $L$-variety.
\end{corollary}

The main thing to do, in order to get the cuspidal
two-orbit varieties of type~II, is to give the list of the
possible $(G,\alpha,\beta)$ where $\alpha$ and $\beta$ are
the positive roots considered previously. 
If $G$ is of rank~$2$, it is done already by Proposition~\ref{rk2}.
So the acting group $G$ will be definitely of rank greater than $2$. 
Sum up the properties of $\alpha$ and $\beta$:

\smallbreak {\rm1)}\enspace
$\alpha\in\Delta$,
$\bigl(w(\lambda),\alpha\bigr)\leq0$ and
$U_\alpha\not\subset G_x$;

\smallbreak {\rm 2)}\enspace
$\lambda-w(\lambda)\in\langle \alpha,\beta\rangle_\mathbb
C$, the $\mathbb C$-vector space spanned by $\alpha$ and
$\beta$;

\smallbreak {\rm 3)}\enspace
$\{\alpha,\beta\}$
basis of $\langle \alpha,\beta\rangle$;

\smallbreak {\rm 4)}\enspace
$\supp\beta\cup\{\alpha\}=\Delta$;

\smallbreak{\rm 5)}\enspace
$w=s_\alpha s_\beta$;

\smallbreak {\rm 6)}\enspace
$(\lambda,\alpha^\vee)=(\lambda,\beta^\vee)$.
\smallbreak
\noindent
It will appear quickly that there are very few roots
satisfying all these conditions mainly because of the two
following statements.

\begin{remark}
Let $\gamma_1$, $\gamma_2$ and $\gamma_3$ be
three roots of $\Phi$ spanning a root system $\Psi$ of
rank~$3$. Suppose $\Psi$ verifies the property:
$$
\eta=n_1\gamma_1+n_2\gamma_2+n_3\gamma_3\in\Psi
\quad
\Longrightarrow
\quad
\eta-n_3\gamma_3\in\Phi
\quad
\mbox{(up to a scalar).}
$$
Then, $\langle\gamma_1,\gamma_2\rangle_\mathbb C \cap
\langle\gamma_3,\eta \rangle_\mathbb C$ is generated by a
root or is equal to $\{0\}$, if $\eta$ is a root of~$\Phi$.
\end{remark}

Let us consider a simple root $\delta$ such that
$(\beta,\delta)>0$ (then $\delta\in\supp\beta$). Set
$\gamma=s_\alpha(\delta)$. 
Then
\begin{equation}
\bigr(w(\lambda),\gamma\bigl)\leq 0.
\label{trick}
\end{equation}

\begin{lemma}
If $\alpha$, $\beta$ and $\gamma$ satisfy the
property given in the remark with $\gamma=\gamma_3$ then
$\bigr(w(\lambda),\gamma\bigl)=0$.
\label{plan}
\end{lemma}

\begin{proof}
If $\bigr(w(\lambda),\gamma\bigl)\neq 0$ then
according to (\ref{trick}),
$\bigr(w(\lambda),\gamma\bigl)<0$. 
And therefore, there exists
$\eta\in\Phi$ such that $\lambda-w(\lambda)\in\langle
\gamma,\eta \rangle_\mathbb C$ 
(argue similarly as we did to get the root $\beta$). 
It implies: $\lambda -w(\lambda
)\in \langle \alpha,\beta \rangle_\mathbb C \cap \langle
\gamma,\eta \rangle_\mathbb C$. But this is impossible
because of the remark and the fact that $T_x^\circ$ is
regular (\emph{see}~\/Proposition~\ref{reg}).
\end{proof}

Start with $G$ classical and suppose: $(\alpha,\delta)=0$.
If $\supp\beta=\Delta$ then $\alpha$, $\beta$ and $\delta$
satisfy the conditions of Lemma~\ref{plan}.
Therefore we get:
$(\lambda,\delta)=(\lambda,\beta)$.
But this equality is incompatible
with $(\lambda,\alpha^\vee)=(\lambda,\beta^\vee)$. So
assume that $\supp\beta=\Delta\setminus\{\alpha\}$.
Then we have the following possibilities for $(G,\alpha,\beta)$:

\smallbreak $\bullet$\enspace
$(A_n,\alpha_i(i=1,n),\tilde{\alpha}-\alpha_i) $;

\smallbreak $\bullet$\enspace
$(G,\alpha_{n},\varepsilon_1-\varepsilon_{n}) $ for
$G=B_{n},C_{n}$;

\smallbreak$\bullet$\enspace
$(G,\alpha_1,\varepsilon_2+\varepsilon_{n}) $ for
$G=B_n,C_n$;

\smallbreak$\bullet$\enspace
$\bigl(G,\alpha_1,\varepsilon_2+\varepsilon_j(2<j<n)\bigr)$
for $G=B_n$, $C_n$, $D_n$;

\smallbreak$\bullet$\enspace
$(D_n,\alpha_{n-1},\varepsilon_1+\varepsilon_{n})$;

\smallbreak $\bullet$\enspace
$(D_n,\alpha_n,\varepsilon_1-\varepsilon_n)$.

\smallbreak
Applying Lemma~\ref{plan}, we end up again with a
contradiction.
Thus necessarily, $(\alpha,\delta)<0$. 
The possible triples $(G,\alpha,\beta)$ are now:

\smallbreak $\bullet$\enspace
$(A_{3},\alpha_2,\tilde{\alpha}) $;

\smallbreak $\bullet$\enspace
$(G,\alpha_2,\varepsilon_1+\varepsilon_{3}) $ for $G=B_n$,
$C_n$ or $D_n$;

\smallbreak $\bullet$\enspace
$(B_{3},\alpha_{3},\tilde{\alpha}) $;

\smallbreak $\bullet$\enspace
$(B_{n},\alpha_1,\varepsilon_2) $;

\smallbreak $\bullet$\enspace
$(C_{n},\alpha_2,\tilde{\alpha}) $;

\smallbreak $\bullet$\enspace
$(C_{n},\alpha_1,2\varepsilon_2) $;

\smallbreak $\bullet$\enspace
$(D_{n},\alpha_1,\tilde{\alpha}) $.

\begin{claim}
In all these cases,
$\bigr(w(\lambda),\gamma\bigl)\neq 0$.
\end{claim}

But because of Lemma~\ref{plan}, we may also have:
$\bigr(w(\lambda),\gamma\bigl)=0$ if $\alpha$, $\beta$ and
$\gamma$ satisfy the good conditions.
So we have to reduce the list to

\smallbreak $\bullet$\enspace
$(A_3,\alpha_2,\tilde{\alpha}) $;

\smallbreak $\bullet$\enspace
$(C_3,\alpha_2,\varepsilon_1+\varepsilon_{3}) $;

\smallbreak $\bullet$\enspace
$(B_3,\alpha_3,\tilde{\alpha}) $;

\smallbreak $\bullet$\enspace
$(C_{n},\alpha_2,\tilde{\alpha}) $; 

\smallbreak $\bullet$\enspace 
$(D_n,\alpha_1,\tilde{\alpha}) $.
\smallbreak

The fourth triple is ruled out just by considering
$\bigr(w(\lambda),\alpha_1\bigl)$.
The other ones give raise to some pairs of Table~2.

For the exceptional case, we proceed similarly and we
obtain the left pairs of Table~2. 
This ends the proof of the main theorem: 
the classification of two-orbit varieties.

\bigbreak

\par\noindent
{\large\textbf{Appendix}}
\bigbreak

To make Table~$3$ and Table~$4$ readable to the reader,
we will need the following notation.

We shall recall, at first, that, in these tables,
the two-orbit varieties are \emph{cuspidal} 
with \emph{connected} generic stabilizers. 
All other two-orbit varieties are obtained
(\emph{see}\/~section~\ref{notation})
either by parabolic induction
or as the quotient of a cuspidal
two-orbit variety by $H/H^\circ$ for $H^\circ$ the corresponding 
connected generic stabilizer given in Tables $1$ and $2$. 

In the first column, we have listed the type of the group $G$ acting
on the two-orbit variety designed in the  second column.
When the action is obvious, we shall not state it precisely.

Once we have the acting group, 
we fix a Borel subgroup $B$ and a maximal torus in it 
and use the standard notation (\emph{see}~\/\cite{BBKI})
to denote $\alpha_i$ the simple roots,
$\omega_i$ the fundamental weight corresponding to $\alpha_i$,
$s_i$ the simple reflection of the Weyl group associated to $\alpha_i$
and $P_i$ the maximal parabolic subgroup attached to $\alpha_i$.

The Grassmaniann of $m$-planes in $\mathbb C^n$ 
is denoted by $\mathop {\rm Gr}(n;m)$. 
For instance,  $\mathop {\rm Gr}(n;1)$ is just $\mathbb P^n$
and  $\mathop {\rm Gr}(n;n-1)$ is $\mathbb P^{n^*}$.

We have a natural action of $SO_n$ on the quadric 
$\mathop{\rm Q}(n)=
\bigl\{[z_0:z_1:\ldots:z_n]\in\mathbb P^n:z_0^2=\sum z_i^2 \bigr\}$
given by $g\cdot[z_0:z']=[z_0:g\cdot z']$, 
$g\in SO_n$ and $z'\in \mathbb C^n$.

As usual, $\Sigma_k$ denotes the Hirzebruch surface: the
$2$-dimensional smooth complete torus embedding corresponding to the
integer $k$.
In particular, when $k=1$, the normalization of $\Sigma_k$ is just
$\mathbb P^2$.


\def\numero{n\raise.82ex\hbox{$\fam0\scriptscriptstyle
o$}~\ignorespaces}

\begin{figure}
\begin{center}
\textbf{Table~3.\enspace Two-orbit varieties of type I}
\par\nobreak\medskip\nobreak
$
\begin{array}{ccll}
A_n
&\triangleright&
\mathbb P^n\times\mathbb P^{n^*}
&
\\
\\
B_n
&\triangleright&
\mathop{\rm Q}(2n+1)
&
\\
\\
B_n
&\triangleright&
\bigl\{(\ell,P)\in\mathbb P^{2n+1}\times\mathop{\rm Gr}(2n+1;n): 
\ell\subset P, \quad\hbox{$P$ totally isotropic}\bigr\}
&
\\
\\
\vcenter{\setbox0\hbox{$B_n$}
        \copy0\medskip\hbox to \wd0{\hss or $C_n$\hss}}
&\triangleright& 
G/P_1\times G/P_n
\\
\\
C_n
&\triangleright&
\mathop{\rm Gr}(2n;2)
&
\\
\\
C_n
&\triangleright&
\bigl\{(\ell,P)\in \mathbb P^{2n}\times\mathop{\rm Gr}(2n;2):
\ell\subset P\bigr\}

\\
\\

F_{4}
&\triangleright&
E_6/P_6
\\
\\
F_4
&\triangleright& 
\vcenter{\setbox0\hbox{
$\displaystyle
\closure (G\cdot x) \subset 
\mathbb P\bigl ( V(\omega_1)\otimes V(\omega_4)\bigr ) $}
\copy0\medskip\hbox to \wd0{with\quad 
$
x=[v_\mu],
\quad \mu=\lambda-\beta, 
\quad \lambda=\omega_1+\omega_4 
$ and 
$\quad  \beta=1111$
\hss}}
\\
\\
G_2 &\triangleright& G\times_{B}\Sigma_3
\\
\\
G_2
&\triangleright& G\times_{B}\Sigma_2
\\
\\
G_2
&\triangleright& G_2\times_{P_1}\mathbb P^2
\quad\quad\hbox{$\mathbb P^2$
compactification of $P_1/H=\mathbb C^2$}
\\
\\
G_2
&\triangleright&
\mathop{\rm Q}(7) \quad\quad\hbox {action induced by $SO_7$'s}

\end{array}
$
\end{center}

\vfill

\begin{center}

\textbf{Table~4.\enspace Two-orbit varieties of type II}
\par\nobreak\medskip\nobreak
$
\begin{array}{ccll}
A_1\times A_1 
&\triangleright& \mathbb P^3 \quad\quad\hbox{compactification of $SL_2$} 
\\
\\
A_2
&\triangleright&
\{(A,z)\in SL_3\times\mathbb P^2: A \hbox{ nilpotent and } Az=0\} 
\\
\\
\vcenter{\setbox0\hbox{$B_2$}
        \copy0\medskip\hbox to \wd0{\hss or $G_2$\hss}}
        &\triangleright&
\closure\bigl( G\cdot[v_\lambda+v_{w(\lambda)}]\bigr)
\subset\mathbb P\bigl(V(\lambda)\bigr) 
\quad\quad\lambda=\omega_1+\omega_2, \quad w=s_1s_2        
\\
\\
B_3
&\triangleright&
\mathop{\rm Q}(8)\quad\quad\hbox {action induced by $SO_8$'s}
\\ 
\\
C_3
&\triangleright&
\Bigl\{ z=\bigl[\sum_{i,j} v_i\wedge v_j\bigr] 
\in\mathbb P(\wedge^2\mathbb C^6): 
\sum_{i,j}\omega(v_i,v_j)=0 \Bigr\}
\\
\\
D_n
&\triangleright& \mathop{\rm Q}(2n)
\\ 
\\
F_4
&\triangleright& 
\vcenter{\setbox0\hbox{
$\displaystyle
\closure (G\cdot x) \subset \mathbb P\bigl ( V(\lambda)\bigr ) $}
\copy0\medskip\hbox to \wd0{with\quad 
$
x=[v_\lambda+v_{w(\lambda)}],
\quad\lambda=\omega_4,
\quad w=s_\alpha s_\beta, 
\quad \alpha=\alpha_4$ and $\beta=1231      
$
\hss}}

\end{array}
$
\end{center}
\end{figure}


\begin{thebibliography}{99}

\bibitem{A1} \textsc{Ahiezer} (D.), \emph{Equivariant
completion of homogeneous algebraic varieties by
homogeneous divisors}, Ann. Global Anal. Geom.~{\bf 1}
(1983), p.~49--78.

\bibitem{A2} \textsc{Ahiezer} (D.), \emph{Dense orbits
with two endpoints}, Izv. Akad. Nauk SSSR Ser. Mat. {\bf
41} (1977), \numero 2, p.~308--324.

\bibitem{Bor} \textsc{Borel} (A.), \emph{Les bouts des
espaces homog\`enes de groupes de Lie}, Ann. of Math. {\bf
58} (1953), p.~443--457.

\bibitem{BBKI} \textsc{Bourbaki} (N.), \emph{Groupes et
alg\`ebres de Lie}, Chapitres 4, 5 et 6, Masson (1981).

\bibitem{BLV} \textsc{Brion} (M.), \textsc{Luna} (D.), \textsc{Vust} (T.),
\emph{Espaces homog\`enes sph\'eriques}, Invent. Math. {\bf 84} (1986),
\numero 3, p.~617--632.

\bibitem{Br1} \textsc{Brion} (M.), \emph{On spherical
varieties of rank one, Group actions and invariant theory},
(Montreal, PQ, 1988), CMS Conf. Proc., {\bf 10}, Amer.
Math.Soc., Providence, RI (1989), p.~31--41.

\bibitem{Br2} \textsc{Brion} (M.), \emph{A note on
two-orbit varieties}, Topology Hawaii (Honolulu, HI, 1990),
World Sci. Publishing, River Edge, NJ (1992), p.~35--40.

\bibitem{BB} \textsc{Bialynicki-Birula} (A.), 
\emph{On fixed point schemes of actions of multiplicative and additive 
groups}, Topology {\bf 12} (1973), p.~99--103.

\bibitem{Cu} \textsc{Cupit-Foutou} (S.),
\emph{Classification des vari\'et\'es \`a deux orbites},
Pr\'epublication de l'Institut de Recherche Math\'ematique Avanc\'ee,
2000/01, www-irma.u-strasbg.fr/irma/publications/2000/00001.shtml.

\bibitem{Feld} \textsc{Feldm\"uller} (D.), \emph{Two-orbit
varieties with smaller orbit of codimension two}, Arch.
Math. (Basel) {\bf 54} (1990), \numero 6, p.~582--593.

\bibitem{H} \textsc{Horrocks} (G.), \emph{Fixed point scheme of additive group actions}, Topology {\bf 8} (1969), p.~233--242.

\bibitem{HS} \textsc{Huckleberry} (A.), \textsc{Snow}
(D.), \emph{Almost-homogeneous K\"ahler manifolds with
hypersurface orbits}, Osaka J. Math. {\bf 19} (1982),
\numero 4, p.~763--786.

\bibitem{Hum} \textsc{Humphreys} (J. E.), \textsc{Linear algebraic groups},
Springer (1981).

\bibitem{LV} \textsc{Luna} (D.), \textsc{Vust} (T.),
\emph{Plongements d'espaces homog\`enes}, Comment. Math.
Helv.~{\bf 58}, (1983), p.~186--245.

\bibitem{Luna} \textsc{Luna} (D.), \emph{Toute vari\'et\'e
magnifique est sph\'erique}, Transform. Groups {\bf 1}
(1996), \numero 3, p.~249--258.

\bibitem{Sumi} \textsc{Sumihiro} (H.), \emph{Equivariant
completion}, J. Math. Kyoto Univ. {\bf 14} (1974),
p.~1--28.

\bibitem{Wass} \textsc{Wassermann} (B.), \emph{Wonderful
varieties of rank two}, Transform. Groups {\bf 1} (1996),
p.~375--403.

\end{thebibliography}
\end{document}